\documentclass[12pt, amsfonts, epsfig]{amsart}
\usepackage{latexsym, amssymb, amscd, amsfonts}

 \newlength{\baseunit}               
 \newcount{\numlines}                
\newcommand{\N}{\mathbb{N}}

\newcommand{\Q}{\mathbb{Q}}

\newcommand{\C}{\mathbb{C}}

        \newfont{\hollow}{msbm10 scaled\magstep1}
        
        \newfont{\Bfmit}{eufm10 scaled\magstep1}

\newcommand{\cO}{{\mathcal O}}
\newcommand{\cL}{{\mathcal L}}

\newcommand{\GL}{\operatorname{GL}}

\newcommand{\Proj}{\operatorname{Proj}}

\newtheorem{thm}{Theorem}[section]

\newtheorem{cor}[thm]{Corollary}

\theoremstyle{definition}

\theoremstyle{remark}

\newcommand{\lremind}[1]{{}}
\newcommand{\bremind}[1]{{}}

\newcommand{\cut}[1]{}

\begin{document}
\pagestyle{plain} \title{{ \large{Factorization Theorem for
Projective Varieties with Finite Quotient Singularities} } }
\author{Yi Hu}

\address{Department of Mathematics, University of Arizona, Tucson, AZ
85721, USA } \email{yhu@math.arizona.edu}

\address{Center for Combinatorics, LPMC,  Nankai University, Tianjin 300071, China}

\begin{abstract} In this paper, we prove that any two birational
projective varieties with finite quotient singularities can be
realized as two geometric GIT quotients of a non-singular
projective variety by a reductive algebraic group. Then, by
applying the theory of Variation of Geometric Invariant Theory
Quotients (\cite{DH}), we show that they are related by a sequence
of GIT wall-crossing flips.
\end{abstract}
\maketitle


{\parskip=12pt 


\section{Statements of Results}

In this paper, we will assume that the ground field is $\C$.

\begin{thm}
\label{thm:GIT}
 Let $\phi: X \rightarrow Y$ be a birational morphism between two  projective
varieties with at worst finite quotient singularities. Then there
is a smooth polarized projective $(\GL_n \times \C^*)$-variety
$(M, \cL)$ such that
\begin{enumerate}
\item $\cL$ is a very ample line bundle and admits two (general)
linearizations $\cL_1$ and $\cL_2$ with $M^{ss}(\cL_1) = M^s
(\cL_1)$ and $M^{ss}(\cL_2) = M^s (\cL_2)$. \item The geometric
quotient $M^s (\cL_1)/(\GL_n \times \C^*)$ is isomorphic to $X$
and the geometric quotient $M^s (\cL_2)/(\GL_n \times \C^*)$ is
isomorphic to $Y$. \item The two linearizations $\cL_1$ and
$\cL_2$ differ only by characters of the $\C^*$-factor, and
$\cL_1$ and $\cL_2$ underly the same linearization of the
$\GL_n$-factor. Let $\underline{\cL}$ be this underlying $\GL_n$-
linearization. Then we have $M^{ss} (\underline{\cL}) = M^s
(\underline{\cL})$.
\end{enumerate}
\end{thm}

As a consequence, we obtain

\begin{thm}
\label{thm:WFT}
 Let $X$ and $Y$ be two birational projective
varieties with at worst finite quotient singularities. Then $Y$
can be obtained from $X$ by a sequence of GIT weighted blowups and
weighted blowdowns.
\end{thm}

The factorization theorem for {\it smooth} projective varieties
was proved by Wlodarczyk and Abramovich-Karu-Matsuki-Wlodarczyka a
few years ago (\cite{AKMW}, \cite{Wlodarczyk00},
\cite{Wlodarczyk03}).  Hu and Keel, in \cite{HK}, gave a short
proof by interpreting it as VGIT wall-crossing flips of
$\C^*$-action. My attention to varieties with finite quotient
singularities was brought out by Yongbin Ruan. The proof here uses
the same idea of \cite{HK} coupled with a key suggestion of Dan
Abramovich which changed the route of my original approach. Only
the first paragraph of \S 2 uses a construction of \cite{HK} which
we reproduce for completeness. The rest is independent. Theorem
\ref{thm:GIT} reinforces the philosophy that began in
\cite{Hu-Keel00}: Birational geometry of $\Q$-factorial projective
varieties is a special case of VGIT.

I thank Yongbin Ruan for asking me about the factorization problem
of projective orbifolds in the summer of 2002 when I visited Hong
Kong University of Science and Technology.  I sincerely thank Dan
Abramovich for suggesting to me to use the results of
Edidin-Hassett-Kresch-Vistoli (\cite{EHKV}) and the results of
Kirwan (\cite{Kirwan85}). I knew the results of \cite{EHKV} and
have had the paper with me since it appeared in the ArXiv, but I
did not realize that it can be applied to this problem until I met
Dan in the Spring of 2004.

\section{Proof of Theorem \ref{thm:GIT}.}

By the construction of \cite{HK} (cf. \S 2 of \cite{Hu-Keel00}) ,
there is a polarized $\C^*$- projective normal variety $(Z, L)$
such that $L$ admits two (general) linearizations $L_1$ and $L_2$
such that
\begin{enumerate}
\item $Z^{ss}(L_1) = Z^s (L_1)$ and $Z^{ss}(L_2) = Z^s (L_2)$.
\item $\C^*$ acts freely on $Z^s(L_1) \cup Z^s (L_2)$.
\item The
geometric quotient $Z^s (L_1)/ \C^*$ is isomorphic to $X$ and the
geometric quotient $Z^s (L_2)/\C^*$ is isomorphic to $Y$.
\end{enumerate}

The construction of $Z$ is short, so we reproduce it here briefly.
Choose an ample cartier divisor $D$ on $Y$. Then there is an
effective divisor $E$ on $X$ whose support is exceptional such
that $\phi^*D = A + E$ with $A$ ample on $X$. Let $C$ be the image
of the injection $\N^2 \to N^1(X)$ given by $(a,b) \to aA + bE$.
The edge generated by $\phi^* D$ divides $C$ into two chambers:
the subcone $C_1$ generated by $A$ and $\phi^* D$, and the subcone
$C_2$ generated by $\phi^* D$ and $E$. The ring $R = \oplus_{(a,b)
\in \N^2} H^0(X, aA + bE)$ is finitely generated and is acted upon
by $(\C^*)^2$ with weights $(a,b)$ on $ H^0(X, aA + bE)$. Let $Z =
\Proj (R)$ with $R$ graded by total degree $(a+b)$. Then a
subtorus $\C^*$ of $(\C^*)^2$ complementary to the diagonal
subgroup $\Delta$ acts naturally on $Z$. The very ample line
bundle $L= \cO_Z(1)$ has two linearizations $L_1$ and $L_2$
descended from two interior integral points in the chambers $C_1$
and $C_2$, respectively. One  verifies (1), (2) by algebra, and
(3) by algebra and the projection formula.

Now, since $\C^*$ acts freely on $Z^s(L_1) \cup Z^s (L_2)$, we
deduce that $Z^s(L_1) \cup Z^s (L_2)$ has at worse finite quotient
singularities. By Corollary 2.20 and Remark 2.11 of \cite{EHKV},
there is a {\it smooth} $\GL_n$- algebraic space $U$ such that the
geometric quotient $\pi: U \rightarrow U/\GL_n$ exists and is
isomorphic to $Z^s(L_1) \cup Z^s (L_2)$ for some $n >0$. Since
$Z^s(L_1) \cup Z^s (L_2)$ is quasi-projective, we see that so is
$U$. In fact, since $Z^s(L_1) \cup Z^s (L_2)$ admits a
$\C^*$-action, all of the above statements can be made
$\C^*$-equivariant. In other words, $U$ admits a $\GL_n \times
\C^*$ action and a very ample line bundle $L_U = \pi^*
(L^k|_{Z^s(L_1) \cup Z^s (L_2)})$ (for some fixed sufficiently
large $k$) with two $(\GL_n \times \C^*)$- linearizations
$L_{U,1}$ and $L_{U,2}$ such that
\begin{enumerate}
\item $U^{ss}(L_{U,1}) = U^s (L_{U,1})$ and $U^{ss}(L_{U,2}) = U^s
(L_{U,2})$. \item The geometric quotient $U^s (L_{U,1})/ (\GL_n
\times \C^*)$ is isomorphic to $X$ and the geometric quotient $U^s
(L_{U,2})/ (\GL_n \times \C^*)$ is isomorphic to $Y$. Moreover,
\item the two linearizations $L_{U,1}$ and $L_{U,2}$ differ only
by characters of the $\C^*$ factor.
\end{enumerate}

Since we assume that $L_U$ is very ample, we have an $(\GL_n
\times \C^*)$- equivariant embedding of $U$ in a projective space
such that the pullback of ${\mathcal O}(1)$ is $L_U$. Let
$\overline{U}$ be the compactification of $U$ which is the closure
of $U$ in the projective space. Let $L_{\overline{U}}$ be the
pullback of ${\mathcal O}(1)$ to $\overline{U}$. This extends
$L_U$ and in fact extends the two linearizations $L_{U,1}$ and
$L_{U,2}$ to $L_{\overline{U},1}$ and $L_{\overline{U},2}$,
respectively,  such that
$$\overline{U}^{ss}(L_{\overline{U},1}) = \overline{U}^s
(L_{\overline{U},1}) = U^{ss}(L_{U,1}) = U^s (L_{U,1})$$ and
$$\overline{U}^{ss}(L_{\overline{U},2}) = \overline{U}^s
(L_{\overline{U},2}) = U^{ss}(L_{U,2}) = U^s (L_{U,2}).$$ It
follows that the geometric quotient $\overline{U}^s
(L_{\overline{U},1}) / (\GL_n \times \C^*)$ is isomorphic to $X$
and the geometric quotient $\overline{U}^s (L_{\overline{U},2}) /
(\GL_n \times \C^*)$ is isomorphic to $Y$.

Resolving the singularities of $\overline{U}$, $(\GL_n \times
\C^*)$-equivariantly, we will obtain a smooth projective variety
$M$. Notice that $\overline{U}^s (L_{\overline{U},1}) \cup
\overline{U}^s (L_{\overline{U},2}) = U^s(L_{U,1}) \cup U^s
(L_{U,2}) \subset U$ is smooth, hence we can arrange the
resolution so that it does not affect this open subset. Let $f: M
\rightarrow \overline{U}$ be the resolution morphism and $Q$ be
any relative ample line bundle over $M$. Then, by the relative GIT
(Theorem 3.11 of  \cite{Hu96}), there is a positive integer $m_0$
such that for any fixed integer $m \ge m_0$, we obtain a very
ample line bundle over $M$, $\cL = f^* L_{\overline{U}}^m \otimes
Q$, with two linearizations $\cL_1$ and $\cL_2$ such that
\begin{enumerate}
\item $M^{ss}(\cL_1) = M^s (\cL_1)= f^{-1}(\overline{U}^s
(L_{\overline{U},1}))$ and $M^{ss}(\cL_2) = M^s (\cL_2) =
f^{-1}(\overline{U}^s (L_{\overline{U},2}))$. \item The geometric
quotient $M^s(\cL_1) / (\GL_n \times \C^*)$ is isomorphic to
$\overline{U}^s (L_{\overline{U},1}) / (\GL_n \times \C^*)$ which
is isomorphic to $X$, and, the geometric quotient $M^s (\cL_2) /
(\GL_n \times \C^*)$ is isomorphic to $\overline{U}^s
(L_{\overline{U},2}) / (\GL_n \times \C^*)$ which is isomorphic to
$Y$.
\end{enumerate}

Finally, we note from the construction that the two linearizations
$\cL_1$ and $\cL_2$ differ only by characters of the
$\C^*$-factor, and $\cL_1$ and $\cL_2$ underly the same
linearization of the $\GL_n$-factor. Let $\underline{\cL}$ be this
underlying $\GL_n$- linearization. It may happen that $M^{ss}
(\underline{\cL}) \ne M^s (\underline{\cL})$. But if this is the
case, we can then apply the method of Kirwan's canonical
desingularization (\cite{Kirwan85}), but we need to blow up
$(\GL_n \times \C^*)$-equivarianly instead of just
$\GL_n$-equivariantly. More precisely, if $M^{ss}
(\underline{\cL}) \ne M^s (\underline{\cL})$, then there exists a
reductive subgroup $R$ of $GL_n$ of dimension at least 1 such that
$$M^{ss}_R (\underline{\cL}):=\{m \in M^{ss}(\underline{\cL}): m \;
\hbox{is fixed by} \; R\}$$ is not empty. Now, because the action
of $\C^*$ and the action of $\GL_n$ commute, using the
Hilbert-Mumford numerical criterion (or by manipulating invariant
sections, or by other direct arguments), we can check that
$$\C^* M^{ss}(\underline{\cL}) = M^{ss}(\underline{\cL}),$$ in
particular,
$$\C^*M^{ss}_R(\underline{\cL}) = M^{ss}_R(\underline{\cL}).$$ Hence, we
have
$$(\GL_n \times \C^*)M^{ss}_R  = \GL_n M^{ss}_R \subset M \setminus M^s (\underline{\cL}).$$
Therefore, we can resolve the singularities of the closure of  the
union of  $\GL_n M^{ss}_R $ in $M$ for all $R$ with the maximal $r
= \dim R$ and blow $M$ up along the proper transform of this
closure. Repeating this process at most $r$ times gives us a
desired nonsingular $(\GL_n \times \C^*)$-variety with
$\GL_n$-semistable locus coincides with the $\GL_n$-stable locus
(see pages 157-158 of \cite{GIT}). Obviously, Kirwan's process
will not affect the open subset $M^{ss}(\cL_1) \cup M^{ss}(\cL_2)=
M^s(\cL_1) \cup M^s(\cL_2) \subset M^s (\underline{\cL})$.  Hence,
this will allow us to assume that $M^{ss} (\underline{\cL}) = M^s
(\underline{\cL})$.

This completes the proof of Theorem \ref{thm:GIT}.

The proof implies the following

\begin{cor}
 Let $\phi: X \rightarrow Y$ be a birational morphism between two  projective
varieties with at worst finite quotient singularities. Then there
is a polarized projective $\C^*$-variety $(\underline{M},
\underline{L})$ with at worst finite quotient singularities such
that $X$ and $Y$ are isomorphic to two geometric GIT quotients of
$(\underline{M}, \underline{L})$ by $\C^*$.
\end{cor}

\section{Proof of Theorem \ref{thm:WFT}}

 Let $\phi: X ---> Y$ be the birational map. By passing to
the (partial) desingularization of the graph of $\phi$, we may
assume that $\phi$ is a birational morphism. This reduces to the
case of Theorem \ref{thm:GIT}.

We will then try to apply the proof of Theorem 4.2.7 of \cite{DH}
(see also \cite{Th}). Unlike the torus case for which Theorem
4.2.7 applies almost automatically, here, because $(\GL_n \times
\C^*)$ involves a non-Abelian group, the validity of  Theorem
4.2.7 must be verified.

From the last section, the two linearizations $\cL_1$ and $\cL_2$
differ only by characters of the $\C^*$-factor, and $\cL_1$ and
$\cL_2$ underly the same linearization of the $\GL_n$-factor. We
denote this common $\GL_n$-linearized line bundle by
$\underline{\cL}$. For any character $\chi$ of the $\C^*$ factor,
let $\cL_\chi$ be the corresponding $(\GL_n \times
\C^*)$-linearization. Note that $\cL_\chi$ also underlies the
$\GL_n$-linearization $\underline{\cL}$. From the constructions of
the compactification $\overline{U}$ and the resolution $M$, we
know that $M^{ss}(\underline{\cL}) = M^s (\underline{\cL})$. In
particular, $\GL_n$ acts with only finite isotropy subgroups on
$M^{ss}(\underline{\cL}) = M^s (\underline{\cL})$.  Now to go from
$\cL_1$ to $\cL_2$,  we will (only) vary the characters of the
$\C^*$-factor, and we will encounter a ``wall'' when a character
$\chi$ gives $M^{ss}(\cL_\chi) \setminus M^s (\cL_\chi) \ne
\emptyset$. In such a case, since $M^{ss}(\cL_\chi) \subset
M^{ss}(\underline{\cL}) = M^s (\underline{\cL})$ which implies
that $\GL_n$ operates on $M^{ss}(\cL_\chi)$ with only finite
isotropy subgroups, the only isotropy subgroups of $(\GL_n \times
\C^*)$ of positive dimensions have to come from the factor $\C^*$,
and hence we conclude that such isotropy subgroups of $(\GL_n
\times \C^*)$ on  $M^{ss}(\cL_\chi)$ have to be one-dimensional
(possibly disconnected) diagonalizable subgroups. This verifies
the condition of Theorem 4.2.7 of \cite{DH} and hence its proof
goes through without changes. (Theorem 4.2.7 of \cite{DH} assumes
that the isotropy subgroup corresponding to a wall is a
one-dimensional (possibly disconnected) diagonalizable group. The
main theorems of \cite{Th} assume that the isotropy subgroup is
$\C^*$ (see his Hypothesis (4.4), page 708).)

\section{GIT on Projective Varieties with Finite Quotient Singularites}

The proof in \S 2 can be modified slightly to imply the following.

\begin{thm}
\label{thm:transform} Assume that a reductive algebraic group $G$
acts on a polarized projective variety $(X, L)$ with at worst
finite quotient singularities. Then there exists a
\underline{smooth} polarized projective variety $(M, \cL)$ which
is acted upon by $(G \times \GL_n)$ for some $n>0$ such that for
any linearization $L_\chi$ on $X$, there is a corresponding
linearization $\cL_\chi$ on $M$ such that $M^{ss}(\cL_\chi)/\!/(G
\times \GL_n)$ is isomorphic to $X^{ss}(L_\chi)/\!/G$. Moreover,
if $X^{ss}(L_\chi) = X^s(L_\chi)$, then $M^{ss}(\cL_\chi) =
M^s(\cL_\chi)$.
\end{thm}

This is to say that all GIT quotients of the {\it singular}
$(X,L)$ ($L$ is fixed) by $G$ can be realized as GIT quotients of
the {\it smooth} $(M,L)$ by $G \times GL_n$. In general, this
realization is a strict inclusion as $(M, \cL)$ may have more GIT
quotients than those coming from $(X,L)$.

When the underlying line bundle $L$ is changed, the
compatification $\overline{U}$ is also changed, so will $M$.
Nevertheless, it is possible to have a similar construction to
include a finitely many different underlying ample line bundles.
However, Theorem \ref{thm:transform} should suffice in most
practical problems because: (1) in most natural quotient and
moduli problems, one only needs to vary linearizations of a fixed
ample line bundle; (2) Variation of the underlying line bundle
often behaves so badly that the condition of Theorem 4.2.7 of
\cite{DH} can not be verified.




\end{document}